\documentclass[11pt]{article}

\usepackage{xyling}
\usepackage{amsmath,amscd}
\usepackage{amsfonts}

\title{An Overview of Hopf Algebras of Trees \protect\\
and Their Actions on Functions}

\author{Robert L. Grossman and Richard G. Larson}

\date{November 15, 2007}

\def\atree#1{ \Tree { \K{$\circ$} \B{d} \\  \K {#1}  } }

\def\aatree#1#2{\Tree { &  \K{$\circ$} \B{dr} \B{dl} \\
\K{#1} && \K{#2} }}

\def\abtree#1#2{\Tree { \K{$\circ$} \B{d} \\
\K {#1} \B{d} \\
\K {#2}}
}

\def\aaatree#1#2#3{\Tree { \K{$\circ$} \B{d} \\
\K{\mbox{${#1}$}} \B{d} \\
\K{ \mbox{${#2}$}} \B{d} \\
\K{\mbox{${#3}$}}
}}

\def\abctree#1#2#3{\Tree { &  \K{$\circ$} \B{dl} \B{d}  \B{dr} \\
\K{\mbox{${#1}$}} & \K{\mbox{${#2}$}} & \K{\mbox{${#3}$}}
}}

\def\abatree#1#2#3{\Tree { &  \K{$\circ$} \B{dl} \B{dr} \\
\K{\mbox{${#1}$}} \B{d} & & \K{\mbox{${#2}$}} \\
\K{\mbox{${#3}$}}
}}

\def\abbtree#1#2#3{\Tree { &  \K{$\circ$} \B{dl} \B{dr} \\
\K{\mbox{${#1}$}} & & \K{\mbox{${#2}$}} \B{d} \\
       & & \K{\mbox{${#3}$}} }}

\def\axbctree#1#2#3{\Tree { &  \K{$\circ$} \B{d} & \\
& \K{\mbox{${#1}$}} \B{dl} \B{dr}  & \\
\K{\mbox{${#2}$}} && \K{\mbox{${#3}$}}
}}

\def\treeplus{\mbox{$\quad$}\raisebox{-5ex}{\mbox{{\large{\bf +}}}}%
   \mbox{$\quad$}}

\def\treeminus{\mbox{$\quad$}\raisebox{-5ex}{\mbox{{\large{\bf -}}}}%
   \mbox{$\quad$}}

\def\treeproduct{\mbox{$\quad$}\raisebox{-5ex}{\mbox{{\large{\bf $\cdot$}}}}%
   \mbox{$\quad$}}

\def\treeequals{\mbox{$\quad$}\raisebox{-5ex}{\mbox{{\large{\bf =}}}}%
   \mbox{$\quad$}}

\def\atreeZ#1{ \Tree { \Kk{1}{$\circ$} \Bk{1}{0}{d} \\  \K {#1}  } }
\def\abtreeZ#1#2{\Tree { \Kk{1}{$\circ$} \Bk{1}{0}{d} \\
\K {#1} \B{d} \\
\K {#2}}
}

\def\bminus{B_{-}}
\def\bplus{B_{+}}

\def\Conn#1#2{\nabla_{#1}{#2}}

\def\coeff#1#2{ a^{\mu_{#2}}_{{#1}} }

\newcommand{\khot}{k\mathcal{T}}
\newcommand{\khotx}[1]{k\mathcal{T}_{#1} }

\newcommand{\Tr}{{\cal T}}
\newcommand{\kT}{k\mathcal{T}}

\newcommand{\Sp}[1]{\mathfrak{S}_{#1}}
\newcommand{\kS}{k\mathfrak{S} }
\newcommand{\kSx}[1]{k\mathfrak{S}_{#1} }
\newcommand{\rmst}[2]{\mathrm{st}(#1,#2)}
\newcommand{\rmstz}[1]{\mathrm{st}(#1)}

\newcommand{\smsh}{\mathbin\#}

\def\Diff#1{\mathrm{Diff}({#1})}
\def\D{\mathcal{D}}
\def\End#1{\mathrm{End}({#1})}

\newtheorem{thm}{Theorem}[section]
\newtheorem{example}[thm]{Example}
\newtheorem{definition}[thm]{Definition}

\begin{document}
\maketitle
\begin{abstract}
We provide an expository account of some of the Hopf algebras
that can be defined using trees, labeled trees, ordered trees
and heap ordered trees.  We also describe some actions of these
Hopf algebras on algebra of functions.
\end{abstract}

\section{Introduction}

There is growing interest in Hopf algebras defined from trees,
permutations, and other combinatorial structures.  In this paper,
we provide an expository account of some of the Hopf algebras
that can be defined using trees, labeled trees, ordered trees
and heap ordered trees.  We also describe some actions of these
Hopf algebras on algebra of functions.

We assume that the reader is familiar with the basic definitions
and properties of Hopf algebras, as covered, for example,
in \cite{moss}.

This paper is based in part on \cite{GL:tree-book}.

\section{Hopf Algebras of Rooted Trees}

\subsection{Products of Trees}

Throughout this paper, $k$ will be a field of characteristic~$0$
such as the real numbers or the complex numbers.

\begin{definition}
By a {\em tree\/} we will mean a finite rooted tree.  Let $\Tr$ be the
set of finite rooted trees, and let $\kT$ be the $k$-vector space
that has $\Tr$ as a basis.
\end{definition}

In this section, we define an algebraic structure on $\kT$
that was introduced in \cite{GL:1989}.
Suppose that $t_1$, $t_2\in\Tr$ are trees.
Let $s_1$, \ldots, $s_r$ be the children of the root of $t_1$
and let $\bminus(\cdot)$ denote the operator that
operates on a tree and removes the root to produce a forest of trees:
\[
\bminus(t_1) = \{ s_1, \ldots, s_r \}.
\]

If $t_2$ has $n+1$ nodes (counting the root), there are
$(n+1)^r$ ways to attach the $r$ subtrees of $t_1$
that have $s_1$, \ldots, $s_r$ as roots to the tree $t_2$.
To do this, make each $s_i$ the child of some node of $t_2$.
We let
\[
\bplus(s_1, \ldots, s_r; t_2).
\]
denote this operation that attaches a forest of trees $s_1$, $\ldots$,
$s_r$ to the nodes of a tree $t_2$ in all possible ways to produce
a set of trees.

The product $t_1t_2$ of the trees $t_1$ and $t_2$ is defined to be the
sum of these $(n+1)^r$ trees.  We summarize this in the following definition:

\begin{definition} The product of the trees $t_1\in\Tr $ and
$t_2\in\Tr$ is defined to be
\[
t_1 t_2 = \sum\bplus\left(\bminus(t_1); t_2 \right)\in \kT.
\]
\end{definition}
Here we sum over the elements of the sets $\bplus$ and $\bminus$.

It is immediate that this definition can be extended by linearity from
$\Tr$ to the vector space $\kT$, and that the trivial tree consisting
only of the root is a right and left unit for this product.  For a
proof that this operation is associative, see \cite{GL:1989}.

Here are some simple examples of tree multiplication:

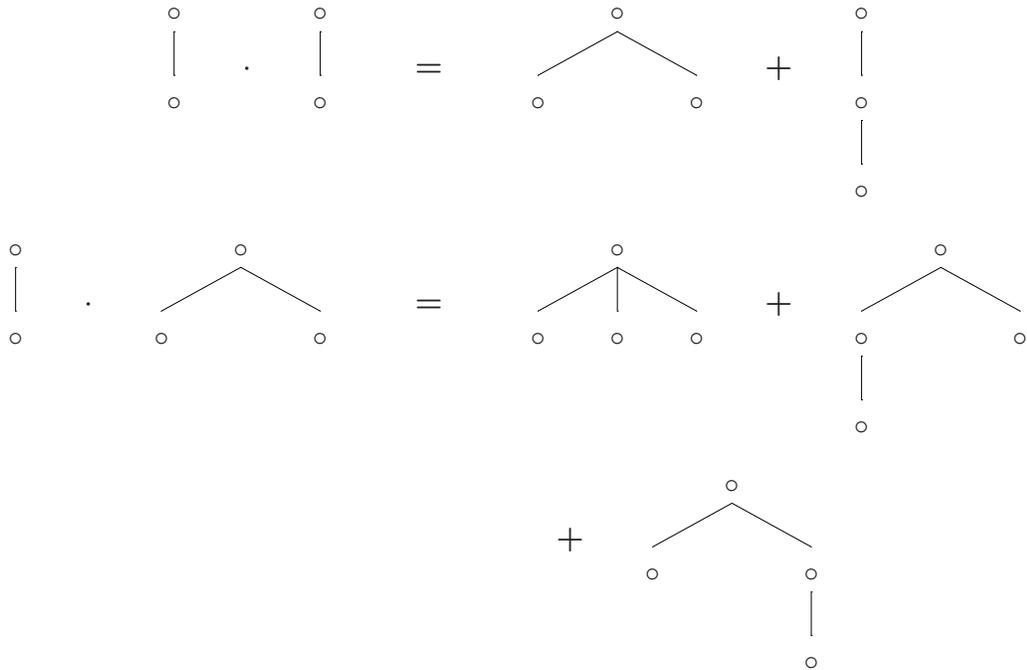
\begin{eqnarray*}
\atree{$\circ$} \treeproduct \atree{$\circ$}
&\treeequals& \aatree{$\circ$}{$\circ$} \treeplus \abtree{$\circ$}{$\circ$} \\
 \atree{$\circ$} \treeproduct \aatree{$\circ$}{$\circ$}
 &\treeequals& \abctree{\circ}{\circ}{\circ} \treeplus \abatree{\circ}{\circ}
  {\circ} \\
  && \quad \treeplus \abbtree{\circ}{\circ}{\circ} \\
\end{eqnarray*}

\subsection{Coproducts of Trees}

We now define a coalgebra structure on the vector space $\kT$:
\begin{eqnarray*}
\Delta(t) &= &\sum_{X\subseteq\bminus(t)}\bplus(X;e)\otimes
   \bplus(\bminus(t)\backslash X;e)                       \\
\epsilon(t) &= &
      \begin{cases}
      1 & \text{if $t$ is the tree whose only node is the root,}  \\
      0 & \text{otherwise.}
      \end{cases}
\end{eqnarray*}
Here, if $X\subseteq Y$ are multisets, $Y\backslash X$ denotes the set
theoretic difference.

It is easily checked that the definition of a
coproduct is satisfied, so that $\kT$ is a cocommutative coalgebra.

This coalgebra and the algebra given in Example~\ref{shuffle-eg}, with
$X$ being the set of trees whose root has exactly one child, are duals
of each other as graded vector spaces, with a tree with $n+1$ nodes
being of degree $n$, and the degree of a monomial $x_{i_1}\cdots
x_{i_p}$ being the sum of the degrees of $x_{i_1}$, \ldots, $x_{i_p}$.
The fact that coassociativity for this coalgebra is trivial gives an
easy proof of the associativity of the algebra in
Example~\ref{shuffle-eg}.

\subsection{Connes--Kreimer Hopf Algebra}

 Another algebra of trees was described in \cite{CK} in 1998 that is
defined as follows:

If ${\cal T}$ is the set of rooted trees, form expressions $a_1 \cdots a_k$,
where $a_i \in {\cal T}$.
This is clearly a commutative algebra of monomials, which we denote $CK$.

The coalgebra structure of $CK$ is defined as follows:
we introduce the notion of {\em cut\/}.
A cut is a set of removed edges.
Given a tree $t$, an {\em admissable cut\/} of $t$ is one with
at most one removed edge on each path from the root to any leaf.
For a cut $C$ we let $R^C(t)$ to be the resulting piece containing the root, and
$P^C(t)$ to be the monomial consisting of the other pieces.
Then define
\[
\Delta(t) = t\otimes 1 + \sum_C P^C(t)\otimes R^C(t).
\]
There is an isomorphism
\[
\chi: \kT \longrightarrow CK^\ast
\]
from the Hopf algebra $\kT$ to the graded dual
of $CK$ satisfying
$${<}\chi(t), a {>} = (\bminus(t), a) = (t, \bplus(a)),$$
where $t\in\kT$ and $a\in CK$.  Here $(\cdot,\cdot)$ is an
inner product defined on $CK$ that depends upon a factor
determined by the symmetry of the arguments.  For the definition of this
inner product and a proof of this fact,
see \cite{Hoffman:2003}.

\subsection{Labeled Trees, Ordered Trees and Heap Ordered Trees}

It is easy to extend the definition of multiplication of trees from
the vector space whose basis is the set of finite rooted trees to the
vector space whose basis is the set of finite {\em labeled} rooted
trees: simply label each node except the root with a label and keep
the labels attached to the nodes when applying the operators $\bminus$
and $\bplus$.  As shown in \cite{GL:1989} this yields a Hopf algebra.

Similarly, the same definition of the product and coproduct applies to
the vector space whose basis is the set of finite {\em ordered trees},
as well as to the vector space whose basis is the set of finite {\em
  ordered, labeled trees}.  Both produce Hopf algebras.  See
\cite{GL:1989} for details.

Later in this paper, we will need the Hopf algebra of heap ordered
trees, which we define following \cite{GL:CA2008}:
\begin{definition}\label{def:hot}
  A standard heap ordered tree on $n+1$ nodes is a finite, rooted tree
  in which all nodes except the root are labeled with the numbers
  $\{1$, $2$, $3$,\ldots, $n\}$ so that:
\begin{enumerate}
\item\label{hot:1} each label $i$ occurs precisely once;
\item\label{hot:2} if a node labeled $i$ has children labeled $j_1$, $\ldots$, $j_k$,
then $i < j_1$,\ldots, $i < j_k$.
\end{enumerate}
\end{definition}

We denote the set of standard heap ordered trees on $n+1$ nodes by
$\mathcal{T}_n$.  Let $\khotx{n}$ be the vector space over the field
$k$ whose basis is the set of trees in $\mathcal{T}_n$, and let
\[
\khot = \bigoplus_{n\ge0}\khotx{n}.
\]

\begin{definition}
A heap ordered tree is a rooted tree in which every node
(including the root) is given a different positive integer label such that
condition~\ref{hot:2}.\ of Definition~\ref{def:hot} is satisfied.
\end{definition}

A heap ordered tree differs from a standard heap ordered tree in that
the root is also labeled, and that the labels can be taken from a
larger set of positive integers.  The number of distinct heap ordered
trees on $n+1$ nodes is $n!$.  Heap ordered trees occur naturally when
studying permutations and differential operators, as we shall see
below.

\subsection{Summary -- Hopf Algebras of Trees}

To summarize, with the product and coproduct as defined above,
the vector spaces whose bases are (i) the set of rooted trees,
(ii) labeled rooted trees, (iii) ordered trees, (iv) labeled ordered
trees, (v) heap ordered trees, and (vi) labeled heap ordered
trees are all Hopf algebras (see~\cite{GL:1989}).

In Section~\ref{section:hot}, we show that there is an isomorphism
from the Hopf algebra of heap ordered trees to a Hopf algebra on
permutations.  In Section~\ref{section:psi1}, we show how Hopf
algebras of rooted labeled trees arise naturally when computing with
derivations on the $k$-algebra of functions from
$k^n$ to $k$.  In
Section~\ref{section:psi2}, we show how Hopf algebras of rooted,
ordered, labeled trees arise naturally when computing with
derivations on the $k$-algebra of $C^{\infty}$ functions on a
$C^{\infty}$ manifold.

\section{Shuffle Algebras}

The definitions and basic properties about shuffle algebras that we
summarize in this section will be needed laer.

A {\em shuffle\/} of the sequences $(i_1,\ldots,i_m)$ and $(j_1,\ldots,j_n)$ is
a permutation $\sigma$ of $(k_1,\ldots,k_{m+n}) = (i_1,\ldots,i_m,j_1,\ldots,j_n)$
satisfying:
if $p<q$, $\sigma(i_p)=k_r$, and $\sigma(i_q)=k_s$, then $r<s$;
if $p<q$, $\sigma(j_p)=k_r$, and $\sigma(j_q)=k_s$, then $r<s$.
In other words, a shuffle of the sequences $(i_1,\ldots,i_m)$ and $(j_1,\ldots,j_n)$
is a permutation of $(i_1,\ldots,i_m,j_1,\ldots,j_n)$ that preserves the order of
$(i_1,\ldots,i_m)$ and of $(j_1,\ldots,j_n)$.

We denote this shuffle by $\sigma(i_1,\ldots, i_m;j_1,\ldots j_n)$.

\begin{definition}[shuffle product]\label{shuffle-eg}
Let $X$ be a set of non-commuting variables,
and let $A$ be the vector space with basis all monomials in the variables
$\{x_1,\ldots,x_n\}$.
Define a product on $A$ as follows.
If $x_{i_1}\cdots x_{i_m}$, $x_{j_1}\cdots x_{j_n}\in A$, define
\[
(x_{i_1}\cdots x_{i_r})(x_{j_1}\cdots x_{j_s})=
\sum_\sigma x_{\sigma(i_1)}\cdots x_{\sigma(i_m)}x_{\sigma(j_1)}\cdots x_{\sigma(j_n)},
\]
where the sum is taken over all shuffles $\sigma(i_1,\ldots, i_r;j_1,\ldots j_s)$.
\end{definition}

The coproduct on $A$ is defined in the same way as for the free non commutative algebra on the set $X$, given by
\[
\Delta(x_i) = 1\otimes x_i + x_i\otimes 1.
\]

\begin{thm}
Fix a field $k$, and let $k{<}x_1, \ldots, x_n{>}$ denote the vector space
over $k$ in the non-commuting variables $x_1$, $\ldots$, $x_n$.  Then
the shuffle product defined above makes $k{<}x_1, \ldots, x_n{>}$ into a
Hopf algebra over $k$.
\end{thm}

\section{Hopf Algebras of Permutations}
\label{section:hot}

Let $\mathfrak{S}_n$ denote the symmetric group on $n$ symbols,
let $\kSx{n}$  denote the vector space over the field $k$ with basis $\mathfrak{S}_n$, and let
\[
\kS = \bigoplus_{n\ge0}\kSx{n}.
\]

In this section, we define a Hopf algebra structure on $\kS$ following
\cite{GL:CA2008}.  We begin with some notation.

Let $(\sigma_1 \sigma_2 \cdots \sigma_k)$ denote the cycle in $\Sp{n}$
which sends $\sigma_1$ to $\sigma_2$, $\sigma_2$ to $\sigma_3$, $\ldots$, and
$\sigma_k$ to $\sigma_1$.
Every permutation is a product of disjoint cycles.
If $(\sigma_1 \sigma_2 \cdots \sigma_k)$ is a cycle, then
there is a string naturally associated with the cycle
that we write $\sigma_1 \sigma_2 \cdots \sigma_k$.

Now let $\sigma=(s_1)\cdots(s_r)\in\Sp{m}$ and
$\tau=(t_1)\cdots(t_\ell)\in\Sp{n}$ be two permutations, each written
as a product of disjoint cycles.  We denote the corresponding strings
by $s_i=m_{i1}\cdots m_{ip_i}$ and $t_j=n_{j1}\cdots n_{jq_j}$
respectively.  We call the elements of $\{n_{11}, \ldots, n_{\ell
  q_\ell}\}$ {\em attachment points} for the cycles $(s_1)$, \ldots,
$(s_r)$ on the permutation $\tau=(t_1)\cdots(t_\ell)$.  We will also
define $\circ$ to be the $(n+1)^{\mathrm th}$ attachment point.

Also, if $\sigma$ is a permutation on $\{1,\ldots,k\}$, let
$\rmst{\sigma}{m}$ be the permutation on $\{m+1,\ldots,m+k\}$ that
sends $m+i$ to $m+\sigma(i)$.

The definition of the product in the bialgebra is simpler if we
introduce the {\em standard order} of a permutation, which is defined as
follows:
\begin{definition}
We say that a permutation $\sigma\in\Sp{m}$ that is expressed
as a product of cycles $\sigma=(s_1)\cdots(s_r)$ is in standard order if
the cycles $(s_i)=(m_{i1}\cdots m_{ip_i})$ are written so that
\begin{enumerate}
\item $m_{i1} < m_{i2}$, $m_{i1} < m_{i3}$, $\ldots$, $m_{i1} < m_{ip_i}$
\item $m_{11} > m_{21} > m_{31}$, $\ldots$, $m_{i-1,1} > m_{i1}$
\end{enumerate}
\end{definition}
In other words, a product of cycles $\sigma=(s_1)\cdots(s_r)$
is written in standard order if
each cycle $(s_i)=(m_{i1}\cdots m_{ip_i})$ starts with its smallest
entry, and if the cycles $(s_1)$, \ldots, $(s_r)$ are ordered so that their
starting entries are decreasing.  A permutation can always be written in
standard order since disjoint cycles commute, and since a single cycle
is invariant under a cyclic permutation of its string.

We now define the {\em heap product\/} of two permutations.

\begin{definition}
We define the heap product $\sigma\smsh\tau$
of $\sigma\in\Sp{m}$ and $\tau\in\Sp{n}$ as follows:
\begin{enumerate}
\item put $\sigma$ and $\tau$ in standard order;
\item replace $\sigma$ by $\rmst{\sigma}{n}$;
\item form terms as follows:
if $(s_i)$ is one of the cycles in $\sigma$, attach the string $s_i$ to any one
of the $n+1$ attachment points;
if the attachment point is one of $n_{11}$,\ldots, $n_{\ell q_\ell}$, say $n_{jk}$,
place the string $s_i=m_{i1}\cdots m_{ip_i}$ to the right of $n_{jk}$; otherwise
(if the attachment point is $\circ$) we multiply
the term we are constructing by $(s_i)$;
\item The product $\sigma\smsh\tau$ is the sum of all the terms
constructed in this way, taken over all the cycles in $\sigma$
and over all attachment points.
\end{enumerate}
\end{definition}
Note that there are $(n+1)^r$ terms in $\sigma\smsh\tau$.

We now describe the coalgebra structure of $\kS$:
Define a function $\rmstz{\pi}$ from permutations to permutations as follows:
let $\pi=(s_1)\cdots(s_p)\in\Sp{n}$ and
let let $L=\{\ell_1,\ldots,\ell_k\}$ be the labels (in order) which occur in the $s_i$.
(If $\pi$ fixes $i$, we include a 1-cycle $(i)$ as a factor in $\pi$.)
The permutation $\rmstz{\pi}$ is the permutation in $\Sp{k}$ gotten by replacing $\ell_j$
with $j$ in $\pi$.
For example, if $\pi=(1 3)(4)(5 7)\in\Sp{7}$, then $\rmstz{\pi}\in\Sp{5}$ equals
$(1 2)(3)(4 5)$.

The coalgebra structure of $\kS$ is defined as follows.
let $\pi=(s_1)\cdots(s_k)\in\Sp{n}$, and let $C=\{(s_1),\ldots,(s_k)\}$.
If $X\subseteq C$ let
$\rho(X) = \rmstz{\prod_{(s_i)\in C}(s_i)}$.
Note that if $\rho(X)\in\Sp{k}$, then $\rho(C\backslash X)\in\Sp{n-k}$.
Define
\begin{eqnarray*}
\Delta(\pi) & = &\sum_{X\subseteq C}\rho(X)\otimes
   \rho(C\backslash X)  \\
\epsilon(\pi) & = &
      \begin{cases}
      1 & \text{if $\pi$ is the identity permutation in $\Sp{0}$,}  \\
      0 & \text{otherwise.}
      \end{cases}
\end{eqnarray*}

The following theorem is proved in \cite{GL:CA2008}.
\begin{thm} With the product and coproduct defined
above, $\kS$ is a Hopf algebra.
\end{thm}

It turns out that this Hopf algebra structure on $\kS$ is isomorphic
to the Hopf algebra structure on heap-ordered trees defined above
\cite{GL:CA2008}.

There are several other different Hopf algebra structures on
permutations that can be defined, including the Malvenuto-Reutenauer
Hopf algebra.  The Malvenuto-Reutenauer is a non-commutative,
non-cocommutative, self-dual, graded Hopf algebra. In contrast, the
Hopf algebra on permutations defined above is non-commutative,
cocommutative, graded Hopf algebra.  For descriptions of these other
Hopf algebras on permutations see \cite{MR:perms},
\cite{HNT:perms-trees}, \cite{AS:MRperms}.

\section{H-Module Algebras for  Labeled Trees  }
\label{section:psi1}

In this and the next section, we describe some actions of Hopf
algebras of trees on $k$-algebras $R$ of functions.  In this section,
we consider the case of $k$-algebras $R$ of functions on a space
$X=k^n$. For example, $R$ can be the space of polynomials, rational
functions or $C^{\infty}$ functions on $X$.  In the next section, we
consider the more general case of $k$-algebras of functions on a
$n$-dimensional algebraic group $X$ or a $n$-dimensional $C^{\infty}$
manifold $X$.

Let $R$ be a commutative $k$-algebra, and let $H$ be a $k$-bialgebra.
The algebra $R$ is a {\em left $H$-module algebra\/} if $R$ is a left
$H$-module for which
\[
h\cdot(ab)=\sum_{(h)} (h_{(1)}\cdot a)(h_{(2)}\cdot b),
\]
where $h\in H$, $\Delta(h)=\sum_{(h)} h_{(1)}\otimes h_{(2)}$,
and $a$, $b\in R$.

Consider the vector space spanned by finite rooted trees, each of
whose nodes except for the root are labeled with one  of the formal
symbols $\{E_1$, \ldots, $E_M\}$.  This generates an algebra using the
product defined above that we denote $\kT(E_1, \ldots, E_k)$ or, if no
confusion is possible, simply by $\kT$.  Note that the same symbol
may label more than one node of the tree.

Fix an $k$-algebra $R$ of functions and
consider the situation where the $E_i$ are not formal symbols,
but are differential operators of the form
\[
E_j = \sum_{\mu=1}^n a_{j}^{\mu}(x) D_{\mu},
\]
with $a_j^{\mu}(x)\in R$.  Here $D_{\mu} = \partial / \partial x_{\mu}$.
The algebra $R$ is typically either the
the algebra of polynomial functions $k[x_1$,
\ldots, $x_n]$, the
algebra of rational functions $k(x_1$, \ldots, $x_n)$,
or the algebra of smooth functions $C^\infty(k^n,k)$.

Let $H$ denote a Hopf algebra of trees labeled with derivations
of $R$.  In this section, we define a natural action of $H$ on $R$ 
that turns $R$ into a $H$-module algebra.

\begin{definition}
We define a map
\[
\psi : \kT(E_1,\ldots,E_M)\longrightarrow {\rm R}
\]
as follows.
Let $t\in\kT(E_1$, \ldots, $E_M)$ have $k+1$ nodes, and let $f\in R$.
Number the non root nodes of $t$ from 1 to $k$.
We define $\psi(t)f$ as follows:
\par\noindent
For the root, form the term
\[
T_0=
\frac{\displaystyle\partial^r f}
{\displaystyle \partial x_{i_l}\cdots\partial x_{i_{l'}}},
\]
where the $r$ children of the root are numbered $l$, \ldots, $l'$.
\par\noindent
For a nonroot node numbered $j$ and labeled with $E_t$, form the term
\[
T_k=
\frac{\displaystyle\partial^r a_t^{j}(x)}
{\displaystyle \partial x_{i_l}\cdots\partial x_{i_{l'}}}.
\]
Then
\[
\psi(t)f =
\sum_{i_1,\ldots,i_k=1}^n T_k T_{k-1}\cdots T_1T_0.
\]
\end{definition}

The following theorem \cite{GL:jsc92} captures an important
property of this $\psi$-map.

\begin{thm}
Fix derivations $E_j = \sum_{\mu} a^{\mu}_j D_{\mu}$ with
  $a^{\mu}_j \in R$, and $D_{\mu} = {{\partial}\over{\partial x_{\mu}}}$.

Let $k{<}E_1,\ldots,E_M{>}$ be the free associative algebra generated
by the formal symbols $E_1$, $\ldots$, $E_M$.
Let $H=\kT(E_1, \ldots, E_M)$ be the Hopf algebra
consisting of rooted trees labeled with the derivations $E_j$.
Let ${\rm Diff}(E_1, \ldots, E_M)$ be the higher order derivations
generated by the derivations $E_1$, $\ldots$, $E_M$.

Then with the $\psi$-map as defined above makes the following
diagram commute:
\[
\begin{array}{rcl}
k{<}E_1,\ldots,E_M{>} & \longrightarrow &\kT(E_1, \ldots, E_M) \\
                      & \searrow        & \;\;\;\downarrow \psi \\
                      &                 & \mathrm{Diff}(E_1, \ldots, E_M)
\end{array}
\]
The right arrow is the map induced by the fact that
$k{<}E_1,\ldots,E_M{>}$ is freely generated by the $E_i$, and $E_i$
maps to the tree with one node other that the root which is labeled
with $E_i$; the diagonal arrow is induced by the fact that the $E_i$
are derivations.
\end{thm}

\begin{example}
An expression such as
\[ p = E_3 E_2 E_1 - E_3 E_1 E_2 - E_2 E_1 E_3 + E_1 E_2 E_3 \]
corresponds to 24 trees, 18 of which cancel (each cancellation saves
$O(n^3)$ differentiations).
The surviving six differentiations are
\begin{eqnarray*}
& & \sum \coeff{3}{3} (D_{\mu_3} \coeff{2}{2}) (D_{\mu_2} \coeff{1}{1}) D_{\mu_1} \\
& & {}- \sum \coeff{3}{3} (D_{\mu_3} \coeff{1}{2}) (D_{\mu_2} \coeff{2}{1}) D_{\mu_1} \\
& & {}-  \sum \coeff{2}{3} (D_{\mu_3} \coeff{1}{2}) (D_{\mu_2} \coeff{3}{1}) D_{\mu_1} \\
& & {}+ \sum \coeff{1}{3} (D_{\mu_3} \coeff{2}{2}) (D_{\mu_2} \coeff{3}{1}) D_{\mu_1} \\
& & {}+ \sum \coeff{3}{3} \coeff{2}{2} (D_{\mu_3} D_{\mu_2} \coeff{1}{1})   D_{\mu_1} \\
& & {}- \sum \coeff{3}{3} \coeff{1}{2} (D_{\mu_3} D_{\mu_2} \coeff{2}{1})   D_{\mu_1}
\end{eqnarray*}
corresponding to the six trees:
\[ \aaatree{E_1}{E_2}{E_3} \treeminus \aaatree{E_2}{E_1}{E_3}
\treeminus \aaatree{E_3}{E_1}{E_2}
\treeplus \aaatree{E_3}{E_2}{E_1}  \]
\[ \treeplus \axbctree{E_1}{E_2}{E_3} \treeminus \axbctree{E_2}{E_1}{E_3} \]
\end{example}

This example illustrates the kinds of cancellations of higher order
derivations that are captured by trees.  It turns out that the Hopf
algebra $H$ of trees labeled with derivations is a natural structure
for computing symbolically with derivations and higher order
derivations.  See \cite{GL:jsc92} and \cite{CG:ISSAC92}.

\section{H-Module Algebras for Labeled, Ordered Trees}
\label{section:psi2}

In this section, we consider what happens when $X$ is not $k^n$, but
instead a $C^{\infty}$ manifold.  In this case, we can define $R$ to
be the $k$-algebra of $k$-valued $C^{\infty}$ functions on the
manifold.

The key observation is that defining natural $H$-module algebras in
this context can be done using a connection and labeled ordered trees.

Let $R$ be a commutative $k$-algebra,
and let $\D$ be a Lie algebra of derivations of $R$.
A {\em connection} \cite{SpivakII} is a map $\D\times\D\rightarrow
\D$ sending $(E,F)\in\D\times\D$ to $\Conn{E}{F}\in\D$
satisfying
\begin{itemize}
\item
$\Conn{E_1+E_2}{F} = \Conn{E_1}{F}+\Conn{E_2}{F}$
\item
$\Conn{E}{(F_1+F_2)} = \Conn{E}{F_1}+\Conn{E}{F_2}$
\item
$\Conn{f\cdot E}{F} = f\cdot\Conn{E}{F}$
\item
$\Conn{E}{(f\cdot F)} = f\cdot\Conn{E}{F}+E(f)F$
\end{itemize}
where $E$, $F\in\D$, $f\in R$.

We use the connection as follows.  If $E,F\in\D$ and $f\in R$, define:
\begin{eqnarray*}
 \atreeZ{E} \!\!\!\! \cdot \left(r\right) & = & E(r)  \\
 \abtreeZ{E}{F} \!\!\!\! \cdot \left(r\right)& = & \Conn{F}{E}(r)
\end{eqnarray*}
and extend using induction, the definition of tree multiplication,
and a consistency requirement for subtrees of larger trees.  We show in
\cite{GL:adv05} that construction defines a $H$-module algebra
structure on $R$.

The following theorem is from \cite{GL:adv05}:
\begin{thm}
Let $R$ be a commutative $k$-algebra, and let $\Conn{E}{F}$ be a
connection on the Lie algebra $\D$ of derivations of $R$.
Then the construction above gives a $\kT(\D)$-module structure on $R$.
This module structure induces a map $\psi:\kT(\D)\longrightarrow\End{R}$.
so that the following diagram commutes:
\[
\begin{array}{rcl}
k{<}E_1,\ldots,E_M{>} & \longrightarrow &\kT(\D) \\
                      & \searrow        & \;\;\;\downarrow \psi \\
                      &                 & \Diff{R}\subset\End{R}
\end{array}
\]
\end{thm}

This theorem can be applied to derive Runge-Kutta numerical algorithms
on groups.  See \cite{CG:JNS93} and \cite{CG:ISSAC92}.

\end{document}